\documentclass[11pt]{article}
\usepackage{amsmath, amssymb, amsthm, graphicx}

\title{On a problem of Erd\H{o}s and Rothschild on edges in triangles}

\author{
Jacob Fox\thanks{Department of Mathematics, MIT,
Cambridge, MA 02139-4307. Email: fox@math.mit.edu. Research
supported by a Simons Fellowship.}
\and
Po-Shen Loh\thanks{Department of Mathematical Sciences, Carnegie Mellon
  University, Pittsburgh, PA 15213. Email: ploh@cmu.edu. Research
  supported by an NSA Young Investigators Grant.}  }
\date{}

\newtheorem{theorem}{Theorem}[section]
\newtheorem*{theorem*}{Theorem}

\newtheorem{lemma}[theorem]{Lemma}
\newtheorem{corollary}[theorem]{Corollary}

\newcommand{\pr}[1]{\mathbb{P}\left[#1\right]}
\newcommand{\E}[1]{\mathbb{E}\left[#1\right]}

\begin{document}
\maketitle

\begin{abstract}
  Erd\H{o}s and Rothschild asked to estimate the maximum number,
  denoted by $h(n,c)$, such that every $n$-vertex graph with at least
  $cn^2$ edges, each of which is contained in at least one triangle,
  must contain an edge that is in at least $h(n,c)$ triangles. In
  particular, Erd\H{o}s asked in 1987 to determine whether for every
  $c>0$ there is $\epsilon>0$ such that $h(n,c)>n^{\epsilon}$ for all
  sufficiently large $n$. We prove that $h(n,c)=n^{O(1/\log \log n)}$
  for every fixed $c<1/4$. This gives a negative answer to the
  question of Erd\H{o}s, and is best possible in terms of the range
  for $c$, as it is known that every $n$-vertex graph with more than
  $n^2/4$ edges contains an edge that is in at least $n/6$ triangles.
\end{abstract}

\section{Introduction}

A \emph{book}\/ of size $h$ in a graph is a collection of $h$
triangles that share a common edge.  The \emph{booksize}\/ of a graph
$G$ is the size of the largest book in $G$. The study of books in
graphs was started by Erd\H{o}s \cite{E62} in 1962, and has since
attracted a great deal of attention in extremal graph theory (see,
e.g., \cite{BN, EFG, EFR, KN}) and graph Ramsey theory (see, e.g.,
\cite{FRS, NR, NR1, NR2, NRS, RS, Su}).

Erd\H{o}s and Rothschild \cite{E87} initiated the study of the
booksize of graphs with the property that every edge is in a
triangle. Let $h(n,c)$ be the largest integer such that every
$n$-vertex graph with at least $cn^2$ edges, each of which is
contained in at least one triangle, must contain an edge that is in at
least $h(n,c)$ triangles. Erd\H{o}s and Rothschild asked to estimate
$h(n,c)$ for fixed $c>0$. This question has received considerable
attention (see, e.g., the Erd\H{o}s problem papers \cite{E87, E88,
  E93}, and the book \cite{CG}).

Using his regularity lemma, Szemer\'edi proved that for every $c>0$,
$h(n,c) \to \infty$ as $n \to \infty$. This fact has a number of
applications to various problems in extremal combinatorics. Ruzsa and
Szemer\'edi \cite{RSz} showed that the statement $h(n,c)>1$ for every
fixed $c>0$ and sufficiently large $n$ implies Roth's theorem: that
every subset of the first $n$ positive integers without a $3$-term
arithmetic progression has size $o(n)$. They also showed that it is
equivalent to the $(6,3)$-theorem: that every $3$-uniform hypergraph
on $n$ vertices in which the union of any $3$ edges contains more than
$6$ vertices has $o(n^2)$ edges.  In the other direction, Alon and
Trotter (see \cite{E93}) proved that for each $c<1/4$ there is $c'>0$
such that $h(n,c)<c'\sqrt{n}$. The condition $c<1/4$ is best possible,
because independent results of Edwards \cite{Ed} and Khad\v ziivanov
and Nikiforov \cite{KN} state that any $n$-vertex graph with more than
$n^2/4$ edges contains an edge in at least $n/6$ triangles. In
particular, this implies for $c>1/4$, we must have $h(n,c) \geq n/6$.

For over two decades, there was no improvement on the $O(\sqrt{n})$
upper bound for any fixed $c<1/4$.  Indeed, Erd\H{o}s even proposed
that perhaps the lower bound should be improved to a power of
$n$. Specifically, in 1987 he asked in \cite{E87} whether there is a
constant $\epsilon>0$ such that $h(n,c)>n^{\epsilon}$ for every fixed
$c>0$ and all sufficiently large $n$. This question was also featured
in the book {\it Erd\H{o}s on Graphs} \cite{CG}. We give a negative
answer to this question. In fact, Theorem \ref{thm:construction} below
implies that $h(n,c)=n^{o(1)}$ for every fixed $c<1/4$. By the above
remark that $h(n,c) \geq n/6$ for $c>1/4$, this gives a best possible
range for $c$ with this bound and shows that a sharp transition occurs
when $c$ is near $1/4$.

\begin{theorem}
  \label{thm:construction}
  For all sufficiently large $n$, there are $n$-vertex graphs with
  $\frac{n^2}{4} \big( 1 - e^{- (\log n)^{1/6}} \big)$ edges, with the
  property that every edge is in a triangle, but no edge is in more
  than $n^{14 / \log \log n}$ triangles.\footnote{All logarithms in this paper are in base $e \approx
2.718$.}
\end{theorem}

The study of $h(n,c)$ with $c$ near $1/4$ began in the problem papers
of Erd\H{o}s \cite{E88, E93}. Let $f$ be such that
$cn^2=n^2/4-f(n)n$. Erd\H{o}s \cite{E88} proved if $f$ is constant,
then $h(n,c)=\Omega(n)$. Bollob\'as and Nikiforov \cite{BN} further
showed that $h(n,c)$ is asymptotically $n/6$ if $f \to 0$. If $f$
tends to infinity with $n$, but not too quickly, so that
$f(n)<n^{2/5}$, they showed that $h(n,c)$ is asymptotically
$\frac{n}{2\sqrt{2f(n)}}$. Note that Theorem \ref{thm:construction}
shows that this behavior cannot continue when $f(n)$ approaches
linearity in $n$. In fact, similar constructions, which we omit, show
that there are positive absolute constants $\alpha,\epsilon$ such that
$h(n,c)=O(n^{1/2-\epsilon})$ where $f(n)=n^{1-\alpha}$. This shows
that the asymptotic behavior of $h(n,c)$ discovered by Bollob\'as and
Nikiforov with $c$ very near $1/4$ already breaks down when $f(n)$ is
some power of $n$ which is less than $1$.

We close the introduction by discussing lower bounds on $h(n,c)$ for
fixed $c>0$.  The fact that $h(n,c)$ tends to infinity follows from
the triangle removal lemma, which is a consequence of Szemer\'edi's
regularity lemma. The triangle removal lemma states that for each
fixed $\epsilon>0$ there is $\delta>0$ such that every graph on $n$
vertices with at most $\delta n^3$ triangles can be made triangle-free
by removing at most $\epsilon n^2$ edges. Suppose $G$ is an $n$-vertex
graph with $cn^2$ edges, each of which is in at least one and at most
$h=h(n,c)$ triangles. The total number of triangles in $G$ is at most
$hcn^2/3$. Thus, if $\delta \geq \frac{hc}{3n}$, then there are
$\epsilon n^2$ edges of $G$ such that every triangle of $G$ contains
at least one of these edges. Since every edge of $G$ is in at least
one triangle, there are at least $cn^2/3$ triangles in $G$, and hence
there is an edge in at least $(cn^2/3)/(\epsilon
n^2)=\frac{c}{3\epsilon} $ triangles. As no edge is in more than $h$
triangles, this implies $\frac{c}{3\epsilon} \leq h$. The regularity
proof gives a bound for $\delta^{-1}$ in the triangle removal lemma
which is a tower of twos of height a power of
$\epsilon^{-1}$. Together with the above bounds on $h(n,c)$, this
implies $h(n,c)$ is at least a power of the iterated logarithm $\log^*
n$. Recently, the first author \cite{F} gave a new proof of the
triangle removal lemma which avoids Szemer\'edi's regularity lemma and
gives a better bound. Namely, in the triangle removal lemma, we can
take $\delta^{-1}$ to be a tower of twos of height logarithmic in
$\epsilon^{-1}$. This gives a lower bound for $h(n,c)$ which is
exponential in $\log^* n$.

\section{Tools}
\label{sec:tools}

The properties of our construction are essentially derived from the
concentration of measure.  Say that a random variable $X(\omega)$ on
an $n$-dimensional product space $\Omega = \prod_{i=1}^n \Omega_i$ is
\emph{$C$-Lipschitz}\/ if changing $\omega$ in any single coordinate
affects the value of $X(\omega)$ by at most $C$.  The Hoeffding-Azuma
inequality (see, e.g., \cite{AS}) provides concentration for these
distributions.

\begin{theorem}[Hoeffding-Azuma Inequality]
  \label{thm:azuma}
  Let $X$ be a $C$-Lipschitz random variable on an $n$-dimensional
  product space.  Then for any $t \geq 0$,
  \begin{displaymath}
    \pr{ | X - \E{X} | > t }
    \leq
    2 \exp\left\{
      -\frac{t^2}{2 C^2 n}
    \right\}.
  \end{displaymath}
\end{theorem}

We also need the following well-known formula for the volume of a
high-dimensional Euclidean ball.  The formula is slightly different
for even and odd dimensions.  Since our analysis is asymptotic in
nature, it suffices to consider only even dimensions (which yield
simpler forms).
\begin{theorem}
  \label{thm:ball-volume}
  For a positive even integer $d$ and a positive real number $r$, the
  volume of $B_r^{(d)}$, the $d$-dimensional Euclidean ball with
  radius $r$, is
  \begin{displaymath}
    \text{Vol}\left(
      B_r^{(d)}
    \right)
    =
    \frac{\pi^{d/2} r^d}{(d/2)!} \,.
  \end{displaymath}
\end{theorem}
\noindent The following weaker estimate turns out to be more
convenient for our analysis. 
\begin{corollary}
  \label{cor:ball-volume}
  For a positive even integer $d$ and a positive real number $r$,
  \begin{displaymath}
    \text{Vol}\left(
      B_r^{(d)}
    \right)
    <
    (2\pi e)^{d/2} \cdot \frac{r^d}{d^{d/2}} \,.
  \end{displaymath}
\end{corollary}

The desired bound in the corollary follows from the standard estimate $d! > \big(
\frac{d}{e} \big)^d$, which is routinely obtained by bounding $\log
(d!) = \sum_{i=1}^d \log i > \int_1^d \log x \, dx$.

\section{Construction}
\label{sec:construction}

We first describe a graph which almost has the desired properties.
Specifically, no edge will be in many triangles, and the number of
edges will be quadratic in the number of vertices, but some edges may
fail to be in triangles.  Throughout this section, we will write $x =
y \pm \delta$ or $x$ is in $y \pm \delta$ to denote $y - \delta \leq x
\leq y + \delta$.

\vspace{2mm}

\noindent \textbf{Pre-Construction.}\, For a positive even integer
$r$, let $d = r^5$, let $n = r^d$, and let $\mu = \frac{r^2 - 1}{6}
\cdot d$.  Consider the tripartite graph with vertex set $A \cup B
\cup C$, where each of $A$ and $B$ are copies of $[r]^d$, and $C =
\{0,1,\ldots,r+1\}^d$.  Vertices $a \in A$ and $b \in B$ are joined by
an edge if and only if (when considered as lattice points in $[r]^d$)
their distance satisfies $\|a-b\|_2^2 = \mu \pm d$.  Similarly,
vertices $b \in B$ and $c \in C$ are adjacent if and only if
$\|b-c\|_2^2 = \frac{\mu}{4} \pm 2d$.  Finally, $c \in C$ and $a \in
A$ are adjacent if and only if $\|c-a\|_2^2 = \frac{\mu}{4} \pm 2d$.

\vspace{2mm}

The following lemma will help us to show that the bipartite graph
between $A$ and $B$ is nearly complete.

\begin{lemma}
  \label{lem:popular-distance}
  Let $r$ and $d$ be given integers, and let $U$ and $V$ be two
  lattice points sampled independently and uniformly at random from
  $[r]^d$.  Define
  \begin{displaymath}
    \mu = \frac{r^2 - 1}{6} \cdot d \,.
  \end{displaymath}
  Then with probability at least $1 - 2 e^{-\frac{d}{2r^4}}$, $\|U -
  V\|_2^2 = \mu \pm d$.
\end{lemma}

\noindent \textbf{Proof.}\, Let $U = (U_1, \ldots, U_r)$ and $V =
(V_1, \ldots, V_r)$.  The squared $L_2$ distance is precisely $\sum_i
(U_i - V_i)^2$, which is a sum of $r$ independent random variables.  A
simple calculation shows that
\begin{displaymath}
  \E{(U_1 - V_1)^2}
  =
  \E{U_1^2} - 2\E{U_1} \E{V_1} + \E{V_1^2}
  = 2 \left( \E{U_1^2} - \E{U_1}^2 \right).
\end{displaymath}
Since $U_1$ is an integer picked uniformly at random from $[r]$, then $\E{U_1} = \frac{r+1}{2}$ while 
\begin{displaymath}
  \E{U_1^2} 
  = 
  \frac{1}{r} \cdot \frac{r(r+1)(2r+1)}{6}
  =
  \frac{(r+1)(2r+1)}{6} \,,
\end{displaymath}
so
\begin{displaymath}
  \E{(U_1 - V_1)^2}
  =
  \frac{r^2 - 1}{6} \,,
\end{displaymath}
and hence
\begin{displaymath}
  \E{\|U - V\|_2^2} = \frac{r^2 - 1}{6} \cdot d \,.
\end{displaymath}
On the other hand, each $(U_i - V_i)^2$ is less than $r^2$, so by
the Hoeffding-Azuma inequality (Theorem \ref{thm:azuma}), the probability that
$\|U-V\|_2^2$ deviates from its expectation by more than $d$ is at
most
\begin{displaymath}
  2 e^{-\frac{d^2}{2r^4 d}}
  =
  2 e^{-\frac{d}{2r^4}} \,,
\end{displaymath}
as claimed.  \hfill $\Box$

\vspace{2mm}

Next, we show that every edge between $A$ and $B$ is in a positive
number of triangles, but not too many.

\begin{lemma}
  \label{lem:AB-in-triangles}
  In the Pre-Construction, the number of edges that join $A$ and $B$
  is at least $(1 - 2e^{-\frac{d}{2r^4}}) n^2$, and every one of those
  edges is contained in between $2^{d-1}$ and $15^d$ triangles.
\end{lemma}

\noindent \textbf{Proof.}\, The first claim is an immediate
consequence of the previous lemma.  We then move to establish a lower
bound on the number of triangles that contain a given edge $ab$.  By
definition, we have $\|a-b\|_2^2 = \mu \pm d$.  Let $m = (m_1, \ldots,
m_d)$ denote the midpoint of $a = (a_1, \ldots, a_d)$ and $b = (b_1,
\ldots, b_d)$ when considered as points in $[r]^d$.  Note that
although $a$ and $b$ have integer coordinates, $m$ may have
half-integer coordinates.  Let $x_i = b_i - a_i$; then $m_i - a_i =
\frac{x_i}{2}$.  For each $i$, if $x_i$ is odd, define $\delta_i =
\frac{1}{2}$, and if $x_i$ is even, define $\delta_i = 1$.  Consider
lattice points $c$ of the form $c_i = m_i + \delta_i \epsilon_i$,
where $\epsilon_i \in \{\pm 1\}$.  All such points still lie in $C$
because $C = \{0, \ldots, r+1\}^d$.  Then,
\begin{align*}
  \|c-a\|_2^2
  &=
  \sum_i \left( \frac{x_i}{2} + \delta_i \epsilon_i \right)^2
  =
  \frac{\|b-a\|_2^2}{4} + \sum_i \delta_i^2 + \sum_i x_i \delta_i \epsilon_i \\
  \|b-c\|_2^2
  &=
  \sum_i \left( \frac{x_i}{2} - \delta_i \epsilon_i \right)^2
  =
  \frac{\|b-a\|_2^2}{4} + \sum_i \delta_i^2 - \sum_i x_i \delta_i \epsilon_i \,.
\end{align*}
Since $\|b-a\|_2^2 = \mu \pm d$ and $\sum_i \delta_i^2 \leq d$, every
choice of $(\epsilon_i)$ satisfying
\begin{displaymath}
  \left| \sum_i x_i \delta_i \epsilon_i \right|
  \leq
  \frac{3}{4}d
\end{displaymath}
will produce a point $c \in C$ which is permissible as the third
vertex of a triangle containing $ab$.  (It would make $\|c-a\|_2^2$
and $\|b-c\|_2^2$ both in $\frac{\mu}{4} \pm 2d$.)  Now consider the
$\epsilon_i$ as independent uniform random variables over $\{\pm 1\}$,
and define the random variable $Z = \sum_i x_i \delta_i \epsilon_i$.
By symmetry, $\E{Z} = 0$, and since $|x_i| \leq r$, changing the
choice of a particular $\epsilon_i$ cannot affect $Z$ by more than
$2r$.  Therefore, the Hoeffding-Azuma inequality (Theorem \ref{thm:azuma}) gives
\begin{displaymath}
  \pr{|Z| > \frac{3}{4} d}
  <
  2 \exp\left\{
    -\frac{ \left(\frac{3}{4} d\right)^2 }{2 (2r)^2 d}
  \right\}
  <
  2 e^{-\frac{d}{15 r^2}} \,,
\end{displaymath}
which implies that the number of valid $c$ is at least
\begin{displaymath}
  \left( 1 - 2 e^{-\frac{d}{15 r^2}} \right) \cdot 2^d 
  > 2^{d-1} \,,
\end{displaymath}
as claimed.

For the upper bound, again assume that we are given $a, b$ such that
$\|a-b\|_2^2 = \mu \pm d$, and let $x_i = b_i - a_i$.  We will bound
the number of half-lattice points $c$ of the form $c_i = a_i +
\frac{x_i}{2} + \frac{w_i}{2}$, where $w_i \in \mathbb{Z}$, which
satisfy $\|c-a\|_2^2 = \frac{\mu}{4} \pm 2d$ and $\|b-c\|_2^2 =
\frac{\mu}{4} \pm 2d$.  For this, observe that
\begin{align*}
  \|c-a\|_2^2
  &=
  \sum_i \left( \frac{x_i}{2} + \frac{w_i}{2} \right)^2
  =
  \frac{\|b-a\|_2^2}{4} + \frac{1}{4} \sum_i w_i^2 
  + \frac{1}{2} \sum_i x_i w_i \\
  \|b-c\|_2^2
  &=
  \sum_i \left( \frac{x_i}{2} - \frac{w_i}{2} \right)^2
  =
  \frac{\|b-a\|_2^2}{4} + \frac{1}{4} \sum_i w_i^2 
  - \frac{1}{2} \sum_i x_i w_i \,,
\end{align*}
so we always have
\begin{displaymath}
  \|c-a\|_2^2 + \|b-c\|_2^2 
  = 
  \frac{\|b-a\|_2^2}{2} + \frac{1}{2} \sum_i w_i^2 \,.
\end{displaymath}
Hence whenever both $\|c-a\|_2^2$ and $\|b-c\|_2^2$ are in $\frac{\mu}{4}
\pm 2d$, we also have $\sum_i w_i^2 \leq 9d$.  It therefore suffices
to bound the number of lattice points in $B_{3 \sqrt{d}}^{(d)}$, the
$d$-dimensional Euclidean ball of radius $3 \sqrt{d}$ centered at the
origin.  Observe that this is at most the volume of
$B_{3.5\sqrt{d}}^{(d)}$, because by placing a unit $d$-dimensional
cube centered at each lattice point in $B_{3 \sqrt{d}}$, we obtain a
non-overlapping collection of unit cubes all contained in the ball of
radius $3 \sqrt{d} + \frac{1}{2} \sqrt{d}$ by the triangle inequality
(the greatest distance from the center of a unit cube to a point on
its boundary is $\frac{1}{2} \sqrt{d}$).

Yet Corollary \ref{cor:ball-volume} bounds the volume of the
$d$-dimensional Euclidean ball of radius $3.5 \sqrt{d}$ by
\begin{displaymath}
  (2 \pi e)^{d/2} 
  \cdot 
  \frac{ \left( 3.5 \sqrt{d} \right)^d }
  { d^{d/2} }
  < 15^d \,,
\end{displaymath}
as claimed.  \hfill $\Box$

\vspace{2mm}

\begin{lemma}
  \label{lem:BC-not-in-many-triangles}
  In the Pre-Construction, every edge joining $B$ and $C$, or joining
  $A$ and $C$, is contained in at most $15^d$ triangles.
\end{lemma}

\noindent \textbf{Proof.}\, Assume that we are given $a, c$ such that
$\|c-a\|_2^2 = \frac{\mu}{4} \pm 2d$, and let $y_i = c_i - a_i$.  We
will bound the number of lattice points $b$ of the form $b_i = a_i + 2
y_i + w_i$, where $w_i \in \mathbb{Z}$, which satisfy $\|b-c\|_2^2 =
\frac{\mu}{4} \pm 2d$ and $\|b-a\|_2^2 = \mu \pm d$.  For this,
observe that
\begin{align*}
  \|b-c\|_2^2
  &=
  \sum_i ( y_i + w_i )^2
  =
  \|c-a\|_2^2 + \sum_i w_i^2 + 2 \sum_i y_i w_i \\
  \|b-a\|_2^2
  &=
  \sum_i ( 2 y_i + w_i )^2
  =
  4 \|c-a\|_2^2 + \sum_i w_i^2 + 4 \sum_i y_i w_i \,,
\end{align*}
and hence
\begin{displaymath}
%  2 \|b-c\|_2^2 + \|b-a\|_2^2
%  =
%  6 \|c-a\|_2^2 + 3 \sum_i w_i^2 \,.
  \|b-a\|_2^2 - 2 \|b-c\|_2^2
  =
  2\|c-a\|_2^2 - \sum_i w_i^2 \,.
\end{displaymath}
Therefore, the only way to have both $\|b-c\|_2^2 = \frac{\mu}{4} \pm
2d$ and $\|b-a\|_2^2 = \mu \pm d$ is to also have $\sum_i w_i^2 \leq
9d$.  By the same computation as in the proof of the previous lemma,
the number of such integral $(w_i)$ is less than $15^d$. Hence, every edge between $A$ and $C$ is 
in at most $15^d$ triangles. By symmetry, every edge between $B$ and $C$ also is in at most $15^d$ 
triangles.  \hfill
$\Box$

\vspace{2mm}

We are now ready to prove the main theorem.

\vspace{2mm}

\noindent \textbf{Proof of Theorem \ref{thm:construction}.}\, Start
with the Pre-Construction for a (sufficiently large) even integer $r$,
with $d = r^{5}$ and $n = r^d=|A|=|B|$.  Note that $n = d^{d/5}$, so
$d = (1 + o(1)) \frac{5 \log n}{\log \log n}$. We will take a random
subgraph by sparsifying $C$.  Let $C' \subset C$ with
$|C'|=2^{-d/2}|C|$ be picked uniformly at random.

Next, consider an edge $ab$ joining $A$ and $B$. By Lemma
\ref{lem:AB-in-triangles}, in the Pre-Construction the edge $ab$ was
in at least $2^{d-1}$ triangles with vertices in $C$. Let $E_{ab}$ be
the event that the edge $ab$ is not in a triangle with a vertex from
$C'$. This happens precisely when none of the $\geq 2^{d-1}$
vertices in $C$ that form a triangle with $ab$ are in $C'$. Hence,
\begin{displaymath}
  \pr{E_{ab}}
  \leq 
  \left.
  \binom{|C|-2^{d-1}}{\frac{|C|}{2^{d/2}}} 
  \middle/ 
  \binom{|C|}{\frac{|C|}{2^{d/2}}} 
  \right.
  \leq \left(
    1 - \frac{2^{d-1}}{|C|}
  \right)^{|C|/2^{d/2}} 
  \leq 
  e^{-2^{\frac{d}{2}-1}} \,,
\end{displaymath}
and the expected number of edges $ab$ for which $E_{ab}$ occurs is at
most
\begin{displaymath}
  |A||B| e^{-2^{\frac{d}{2}-1}}
  =
  n^2 e^{-2^{\frac{d}{2}-1}} \,.
\end{displaymath}
Fix a choice of $C'$ with at most $n^2 e^{-2^{\frac{d}{2}-1}}$ edges
$ab$ satisfying $E_{ab}$. Consider the subgraph induced by $A \cup B
\cup C'$.  The total number of vertices in the graph is only
\begin{equation}
  \label{eq:N}
  N 
  = 2n+|C'|=\left(2+2^{-d/2}\left(\frac{r+2}{r}\right)^d\right)n
  <
  \left(
    2 + 2^{- \frac{d}{2}} \cdot e^{\frac{2d}{r}}
  \right) n
  <
  \left(
    2 + 2^{- \frac{d}{3}}
  \right) n \,.
\end{equation}

Unfortunately, now some edges are no longer in triangles.  We resolve
this by deleting all such edges.  By Lemma \ref{lem:AB-in-triangles}
with $r = d^{1/5}$, the number of edges between $A$ and $B$ was
originally at least $n^2 \big( 1 - 2 e^{- \frac{1}{2} d^{1/5}} \big)$,
so since we chose $C'$ such that at most $n^2 e^{-2^{\frac{d}{2}-1}}$
edges $ab$ are not in triangles, the number of remaining edges between
$A$ and $B$ after deleting those not in triangles is still at least
$n^2 \big( 1 - 3 e^{- \frac{1}{2} d^{1/5}} \big)$.  Therefore, by
\eqref{eq:N}, the number of remaining edges between $A$ and $B$ is at
least
\begin{displaymath}
  \frac{N^2}
  { \left( 2 + 2^{- \frac{d}{3}} \right)^2 }
  \left( 1 - 3 e^{- \frac{1}{2} d^{1/5}} \right)
  >
  \frac{N^2}{4} \left(
    1 - 4 e^{- \frac{1}{2} d^{1/5}}
  \right)
  >
  \frac{N^2}{4} \left(
    1 - e^{- (\log N)^{1/6}}
  \right) \,,
\end{displaymath}
so the remaining graph has the claimed total number of edges.
Finally, note that our deletions cannot create any new triangles, so
by Lemmas \ref{lem:AB-in-triangles} and
\ref{lem:BC-not-in-many-triangles}, every edge is still in at most
\begin{displaymath}
  15^d
  =
  15^{(1 + o(1)) \frac{5 \log n}{\log \log n}}
  <
  N^{14 / \log \log N}
\end{displaymath}
triangles, completing our proof.  \hfill $\Box$

\vspace{2mm}

\noindent \textbf{Remark 1.}\, The use of randomness to pick $C'$ in
the above construction is not necessary. Indeed, the construction can
be made explicit by instead picking $C'$ greedily so that each new
vertex added to $C'$ (locally) maximizes the number of edges between
$A$ and $B$ that are in triangles with vertices from $C'$.

\vspace{2mm}

\noindent \textbf{Remark 2.}\, After publicizing this result, the
authors received the following nice observation from Noga Alon.  The
objective of sparsifying $C$ to $C'$ was to raise the edge density to
approach $1/4$.  A simpler way to increase the density is to leave $C$
alone, and instead replace each vertex of $A \cup B$ with exactly
$2^d$ copies of itself, joining two copies of (different) vertices by
an edge if their original vertices were initially adjacent, and
joining a copy of a vertex in $A \cup B$ to an uncopied vertex $c \in
C$ if the corresponding original vertex of $A \cup B$ was adjacent to
$c$.  This avoids our final probabilistic arguments altogether, and
allows for the further simplification that in the Pre-Construction,
all of $A, B, C$ can be taken to be $[r]^d$.  Then, it suffices to
replace the lower bound in Lemma \ref{lem:AB-in-triangles} with the
observation that for any edge $ab$ between $A$ and $B$, the
integer-rounded midpoint produces at least one point $c$ which
completes $ab$ to a triangle.

\section*{Acknowledgment}

We would like to thank Noga Alon for helpful discussions.

\end{document}